\documentclass[reqno,12pt]{amsart}

\numberwithin{equation}{section}

\usepackage{amsmath}

\usepackage{esint} 

\usepackage{amsthm}

\usepackage{verbatim}

\usepackage{calligra}

\usepackage{epsfig}

\usepackage{psfrag}

\usepackage{graphicx}

\usepackage{graphpap,latexsym,epsf}

\usepackage{color}

\usepackage{amssymb,mathrsfs,enumerate}

\usepackage{endnotes}

\usepackage{calligra}

\usepackage{mathtools}

\usepackage[colorinlistoftodos]{todonotes} 

\definecolor{citegreen}{rgb}{0,0.6,0}

\definecolor{refred}{rgb}{0.8,0,0}

\usepackage[colorlinks, citecolor=citegreen, linkcolor=refred]{hyperref}

\usepackage{enumitem}

\newcommand{\R}{\mathbb{R}}

\def\RRR{{\mathrm R}}

\newcommand{\pa}{\partial}

\newcommand{\ep}{\varepsilon}

\newcommand{\Ric}{{\rm Ric}}

\newcommand{\De}{\Delta}

\newcommand{\na}{\nabla}

\mathchardef\emptyset="001F

\definecolor{vgreen}{rgb}{0.1,0.5,0.2}

\definecolor{viola}{RGB}{85,26,139}

\newtheorem{theorem}{Theorem}[section]

\newtheorem{remark}[theorem]{Remark}

\newtheorem{ackn}{Acknowledgements
\hspace{-.4cm}
}

\definecolor{byzantium}{rgb}{0.44, 0.16, 0.39}

\definecolor{amber}{rgb}{1.0, 0.75, 0.0}

\definecolor{darkmagenta}{rgb}{0.55, 0.0, 0.55}

\definecolor{fuzzywuzzy}{rgb}{0.8, 0.4, 0.4}

\definecolor{brown}{rgb}{0.2, 0.08, 0.08}

\definecolor{arancio}{rgb}{1.0, 0.13, 0.0}

\begin{document}
\title[Estimates 
for
the $p$-Green function near the pole]{Estimates 
for
the $p$-Green function near the pole}
\date{}
\author[V.~Agostiniani]{Virginia Agostiniani}
\address{V.~Agostiniani, Universit\`a degli Studi di Trento,
via Sommarive 14, 38123 Povo (TN), Italy}
\email{virginia.agostiniani@unitn.it}
\author[S.~Borghini]{Stefano Borghini}
\address{S.~Borghini, Universit\`a degli Studi di Napoli Federico II,
Via Cintia, Monte S. Angelo, 80126 Napoli (NA), Italy}
\email{stefano.borghini@unina.it}

\begin{abstract}
We study the asymptotic expansion of $p$-Green functions and their derivatives near the pole. In the Euclidean setting, we strengthen a theorem from~\cite{AnEs} by establishing improved integrability properties and estimates for the first and second derivatives. In the more general Riemannian setting, we derive a refined asymptotic expansion of the $p$-Green function near its pole, thereby improving the corresponding result proved in~\cite{Ben_Mar_Rig_Set_Xu}.
\end{abstract}

\maketitle

\bigskip




\section{Introduction}

In this paper, we study the asymptotic expansion of the $p$-Green function near its pole, on a Riemannian manifold of dimension $n\geq 2$. In particular, we aim to describe as precisely as possible the sense in which any $p$-Green function is asymptotic to the model rotationally symmetric Euclidean $p$-Green function, given by
\begin{equation}
\label{eq:Gamma}
\Gamma(x)\,=\,\frac{p-1}{n-p}\,|x|^{-\frac{n-p}{p-1}}\,.
\end{equation}

When the ambient manifold is the Euclidean space, it has been proved in~\cite{Serrin} that any $p$-Green function is asymptotic to $\Gamma$ at the pole. This result has been strenghtened in~\cite{KiVe} and then recently in~\cite{AnEs}, where it has been proved that there exists $\alpha\in(0,1)$ such that the error between a $p$-Green function and $\Gamma$ is $C^{0,\alpha}$ at the pole. We obtain the following improvement on the control of the error near the pole:

\begin{theorem}
\label{thm:Euclidean_1}
Let $\Omega\subset\R^n$ be a bounded domain containing the origin $o$, and let $p\in(1,n)$. Let $u$ be a function satisfying
\[
\De_p^{(g_{\R^n})} u\,=\,-|\mathbb{S}^{n-1}|\,\delta_o \quad \hbox{in }\Omega
\]
and let $\Gamma$ be given by~\eqref{eq:Gamma}.
Then there exists $\alpha\in(0,1)$ such that the function $e=u-\Gamma$ is $C^{0,\alpha}$ at the origin and satisfies
\[
e(x)-e(o)=
\mathcal{O}_2(|x|^\alpha)\,.
\]
\end{theorem}

Referring the reader to Remark \ref{rmk:O} for a definition of $\mathcal O_2$, we note 
in particular that the gradient of the error satisfies $|\na e|=\mathcal{O}(|x|^{\alpha-1})$. It follows immediately that $\na e$ belongs to $L^q$ for all $q\leq n/(1-\alpha)$. In particular, it belongs to $L^n$. It is interesting to notice that this integrability property does not depend on $p$. Furthermore, it improves on~\cite{AnEs}, where it is proved that $\na e$ belongs to $L^{\bar q}$, with $\bar q=n(p-1)/(n-1)$, which is less than $n$ since $p< n$.

\medskip

For general Riemannian manifolds, the error estimates in~\cite{AnEs, KiVe} do not always hold. In fact, as pointed out in~\cite[Remark 2.5]{Mar_Rig_Set}, on the round sphere and on the hyperbolic space the difference between a $p$-Green function and $\Gamma$ can go to $+\infty$ when $p<(n+2)/3$, see the computation in Appendix~\ref{appendix}.
On the other hand, we are able to obtain the following:

\begin{theorem}
\label{thm:Riemannian_1}
Let $(M,g)$ be a $n$-dimensional Riemannian manifold and let $p\in(1,n)$. Let $o\in M$ and let $u$ be a function satisfying
\[
\De^{(g)}_p u\,=\,-|\mathbb{S}^{n-1}|\,\delta_o
\]
in an open neighborhood of $o$.
Let $\{x^i\}$ be normal coordinates centered at $o$ and let $\Gamma$ be given by~\eqref{eq:Gamma}.
Then the error $e=u-\Gamma$ satisfies
\[
e(x)\,=\,
\begin{cases}
\mathcal{O}_2(1) & \hbox{if }
\frac{n+2}{3}<p<n
\\
\mathcal{O}_2(\log(1/|x|)) & \hbox{if }p=\frac{n+2}{3}
\\
\mathcal{O}_2(|x|^{-\frac{n+2-3p}{p-1}}) & \hbox{if }1<p<\frac{n+2}{3}
\end{cases}
\]
as $|x|\to 0$.
\end{theorem}

\begin{remark}
\label{rmk:O}
With the notation $e(x)=\mathcal{O}_2(f(|x|))$ we mean that $e(x)=\mathcal{O}(f(|x|))$, that the first derivatives of $e(x)$ are $\mathcal{O}(f(|x|)/|x|)$ and the second derivatives are $\mathcal{O}(f(|x|)/|x|^2)$ as $|x|\to 0$. 
\end{remark}

Theorem~\ref{thm:Riemannian_1} improves on~\cite[Theorem~2.1]{Ben_Mar_Rig_Set_Xu}, which states that
\[
e(x)\,=\,o_2(|x|^{-\frac{n-p}{p-1}})\,,
\]
where the notation $o_2$ is defined analogously to $\mathcal{O}_2$.

Notice that Theorem~\ref{thm:Riemannian_1} shows in particular that the error is bounded when $p\in((n+2)/3,n)$. It would be interesting to investigate whether the error $e$ is H\"older continuous for this range of values of $p$, in analogy with the Euclidean case.

The asymptotic estimates obtained in Theorem~\ref{thm:Riemannian_1} are sharp for all $p\neq (n+2)/3$, as one can check looking at the behaviour of $e(x)$ on space forms (see again Appendix~\ref{appendix}). In the case $p=(n+2)/3$, we expect, again comparing with space forms, that it should be possible to improve on the asymptotic estimates for the first and second derivatives. Namely, for $p=(n+2)/3$ we expect the first derivative of $e$ to decay as $1/|x|$ and the second derivative to decay as $1/|x|^2$. 

The estimates in Theorem~\ref{thm:Riemannian_1} are crucial in~\cite{AgBoMa}, where they are used to justify the physical relevance of suitable mass-like quantities, by showing that they satisfy a small sphere limit property.

Finally, we remark that in our main theorems we do not cover the case $p=n$ (that is instead considered in~\cite{AnEs,Ben_Mar_Rig_Set_Xu}), as the strategy employed in Section~\ref{sec:estimates} appears to be hard to adapt to this case. Nevertheless, in Section~\ref{sec:order_zero} we prove that, in the setting of Theorem~\ref{thm:Riemannian_1}, if $p=n$ it holds $e(x)=\mathcal{O}(1)$, improving on the estimate in~\cite{Ben_Mar_Rig_Set_Xu} (however we do not discuss estimates for the derivatives of $e$).

\section{Zeroth-order estimate of the error}
\label{sec:order_zero}

In the Euclidean setting, it is known from~\cite{KiVe} that the error term is bounded near the pole.
For general Riemannian manifolds, it is known that $\lim_{|x|\to 0}u(x)/\Gamma(x)= 1$ (see~\cite{Ben_Mar_Rig_Set_Xu,Serrin}), however the error $e$ is not necessarily bounded. On the other hand, in this section we show that the error term $e(x)=u(x)-\Gamma(x)$ satisfies
\begin{equation}
\label{eq:assumption_error}
e(x)\,=\,
\begin{dcases}
\mathcal{O}(1) & \hbox{if } \frac{n+2}{3}<p\leq n\,,
\\
\mathcal{O}(\log 1/|x|)  & \hbox{if } p=\frac{n+2}{3}\,,
\\
\mathcal{O}(|x|^{-\frac{n+2-3p}{p-1}}) & \hbox{if }1<p<\frac{n+2}{3}\,.
\end{dcases}
\end{equation}
We remark that in this section we consider the case $p=n$ as well, in which case we denote as $\Gamma(x)$ the function $\log(1/|x|)$. We also emphasize that we will actually prove something more than~\eqref{eq:assumption_error}, namely we will obtain explicit bounds on the error $e$ in terms of the Ricci curvature of the metric at the pole.
The proof proposed in this section is an adaptation of the proof of the boundedness of the error in the Euclidean setting done in~\cite{KiVe}.
The same circle of ideas has been employed in the proof of~\cite[Theorem~2.1]{Ben_Mar_Rig_Set_Xu}, where the authors obtain estimates on the error $e$ in terms of the sectional curvature of the metric at the pole. 

To prove~\eqref{eq:assumption_error}, we start by considering the function
\[
\Gamma_\kappa(x)\,=\,
\begin{dcases}
\frac{p-1}{n-p}|x|^{-\frac{n-p}{p-1}}+\frac{\kappa}{6(n+2-3p)}|x|^{-\frac{n+2-3p}{p-1}} & \hbox{if }p\in(1,n)\setminus\left\{\frac{n+2}{3}\right\}\,,
\\
\frac{1}{2}|x|^{-2}+\frac{\kappa}{2(n-1)}\log\frac{1}{|x|} & \hbox{if }p=\frac{n+2}{3}\,,
\\
\log\frac{1}{|x|}-\frac{\kappa}{12(n-1)}|x|^2  & \hbox{if }p=n\,.
\end{dcases}
\]
These formulas correspond to the truncated expansions near the origin of the rotationally symmetric $p$-harmonic function on the space form with sectional curvature equal to $\kappa/(n-1)$, see Appendix~\ref{appendix}.

Let us recall that the metric coefficients in normal coordinates
can be expanded as
\begin{equation}
\label{metric_expansion}
g_{ik}(x)\,=\,\delta_{ik}-\frac{1}{3}\RRR_{ijkl}(o)x^j x^l+\mathcal O_1(|x|^3)\,,
\end{equation}
where we have used Einstein notation (i.e, summation over repeated indices is implied) and the following convention for the Riemann tensor:
\[
\RRR(X,Y)Z\,=\,
\nabla_Y\nabla _XZ-\nabla_X\nabla_YZ-\nabla_{[Y,X]}Z\,,
\]
$X,Y,Z$ being vector fields
(see, e.g.,  \cite[Osservazione 5.1.9 and Proposizione 4.6.23]{Man_Maz_Oro_book}
for a proof of the above expansion).
Formula \eqref{metric_expansion} 
easily yields
\[
g^{ik}(x)\,=\,\delta_{ik}+\frac{1}{3}\RRR_{ijkl}(o)x^j x^l+\mathcal{O}_1(|x|^3)\,,
\qquad
\sqrt{g}(x)\,=\,1-\frac{1}{6}\RRR_{jl}(o)x^j x^l+\mathcal{O}_1(|x|^3)\,,
\]
and in turn
\[
\frac{\pa g^{ik}}{\pa x^i}(x)\,=\,-\frac{1}{3}\RRR_{jk}(o)x^j+\mathcal{O}(|x|^2)\,,\quad \frac{\pa\sqrt{g}}{\pa x^i}(x)\,=\,-\frac{1}{3}\RRR_{ij}(o)x^j+\mathcal{O}(|x|^2)\,.
\]

Although the formula for $\Gamma_\kappa$ changes depending on $p$, one can check that, for all $p\in(1,n]$, it holds
\begin{equation}
\label{eq:dGamma}
\frac{\pa\Gamma_\kappa}{\pa x^i}(x)\,=\,-|x|^{-\frac{p+n-2}{p-1}}\left[1+\frac{\kappa}{6(p-1)}|x|^2\right] x^i\,.
\end{equation}
Since all the following formulas depends only on the derivatives of $\Gamma_\kappa$, we can discuss all values of $p$ simultaneously. For all $p\in(1,n]$, we compute
\begin{align*}
|d\Gamma_\kappa|_g^2(x)&=|x|^{-2\frac{n-1}{p-1}}\left[1+\frac{\kappa}{3(p-1)}|x|^2+\mathcal{O}(|x|^3)\right]\,,
\\
|d\Gamma_\kappa|_g^{p-2}\frac{\pa\Gamma_\kappa}{\pa x^i}(x)&=-|x|^{-n}\left[1+\frac{\kappa}{6}|x|^2+\mathcal{O}(|x|^3)\right]x^i\,,
\\
\frac{\pa}{\pa x^i}\!\left[|d\Gamma_\kappa|_g^{p-2}\frac{\pa\Gamma_\kappa}{\pa x^k}\right]\!(x)&=|x|^{-n}\!\left[\left(-\delta_{ik}+n\frac{x^i x^k}{|x|^2}\right)\!+\frac{\kappa}{6}\left(-\delta_{ik}+(n-2)\frac{x^i x^k}{|x|^2}\right)\!|x|^2\!+\mathcal{O}(|x|^3)\right],
\end{align*}
In particular, we obtain
\begin{align*}
g^{ik}\frac{\pa}{\pa x^i}\left[|d\Gamma_\kappa|_g^{p-2}\frac{\pa\Gamma_\kappa}{\pa x^k}\right](x)\,&=\,|x|^{2-n}\left[-\frac{\kappa}{3}-\frac{1}{3}\RRR_{jl}(o)\frac{x^j x^l}{|x|^2}+\mathcal{O}(|x|)\right]\,,
\\
\frac{\pa g^{ik}}{\pa x^i}|d\Gamma_\kappa|_g^{p-2}\frac{\pa\Gamma_\kappa}{\pa x^k}(x)\,&=\,\frac{1}{3}|x|^{2-n}\RRR_{jl}(o)\frac{x^j x^l}{|x|^2}\left[1+\mathcal{O}(|x|)\right]\,,
\\
\frac{1}{\sqrt{g}}\frac{\pa \sqrt{g}}{\pa x^i} g^{ik}|d\Gamma_\kappa|_g^{p-2}\frac{\pa\Gamma_\kappa}{\pa x^k}(x)\,&=\,\frac{1}{3}|x|^{2-n}\RRR_{jl}(o)\frac{x^j x^l}{|x|^2}\left[1+\mathcal{O}(|x|)\right]\,.
\end{align*}
Finally, recall that the $p$-Laplacian can be computed as
\begin{align*}
\De_p^{(g)}\Gamma_\kappa\,&=\,\frac{1}{\sqrt{g}}\frac{\pa}{\pa x^i}\left[\sqrt{g}\, g^{ij}\,|d\Gamma_\kappa|_{g}^{p-2}\frac{\pa\Gamma_\kappa}{\pa x^j}\right]
\\
&=\,\frac{1}{\sqrt{g}}\frac{\pa\sqrt{g}}{\pa x^i} g^{ij}|d\Gamma_\kappa|_{g}^{p-2}\frac{\pa\Gamma_\kappa}{\pa x^j}+
\frac{\pa g^{ij}}{\pa x^i}|d\Gamma_\kappa|_{g}^{p-2}\frac{\pa\Gamma_\kappa}{\pa x^j}+
g^{ij}\frac{\pa}{\pa x^i}\left[|d\Gamma_\kappa|_{g_s}^{p-2}\frac{\pa\Gamma_\kappa}{\pa x^j}\right]\,.
\end{align*}
Putting together the previous estimates, we conclude
\begin{equation*}
\De_p^{(g)}\Gamma_\kappa\,
=\,\left(\frac{\RRR_{jl} x^j x^l}{|x|^2}-\kappa\right)\frac{1}{3|x|^{n-2}}+\mathcal{O}(|x|^{3-n})\,.
\end{equation*}
Thus, if $a,b\in\R$ are such that
\begin{equation}
\label{eq:Ricci_bound}
a\, g(o)\,\leq\, \Ric(o) \,\leq\, b\, g(o)\,,
\end{equation}
then we have
\begin{align*}
\De_p^{(g)}\Gamma_{a-\ep}(x)\,&\geq\,\frac{\ep}{3|x|^{n-2}}\left[1+\mathcal{O}(|x|)\right]\,,
\\
\De_p^{(g)}\Gamma_{b+\ep}(x)\,&\leq\,-\frac{\ep}{3|x|^{n-2}}\left[1+\mathcal{O}(|x|)\right]\,.
\end{align*}
Thus, there exists $r_\ep>0$ such that on the whole $B_{r_\ep}(o)$ it holds
\[
\De_p^{(g)}\Gamma_{a-\ep}(x)\,\geq\,0\,,\qquad
\De_p^{(g)}\Gamma_{b+\ep}(x)\,\leq\,0\,.
\]
Now let $0<\delta<1$ and consider the two functions
\begin{align*}
v_\delta^-(x)\,&=\,(1-\delta)\Gamma_{a-\ep}(x)-(1-\delta)\Gamma_{a-\ep}(r_\ep)+\inf_{\{|x|=r_\ep\}}u\,,
\\
v_\delta^+(x)\,&=\,(1+\delta)\Gamma_{b+\ep}(x)-(1+\delta)\Gamma_{b+\ep}(r_\ep)+\sup_{\{|x|=r_\ep\}}u\,.
\end{align*}
where we have used the slightly imprecise notation $\Gamma_\kappa(r_\ep)$ to denote the constant value of $\Gamma_\kappa(x)$ on the points with $|x|=r_\ep$.
By construction, we have $v_\delta^-\leq u\leq v_\delta^+$ on $\pa B_{r_\ep}(o)=\{|x|=r_\ep\}$. Furthermore, since $\lim_{|x|\to 0}u(x)/\Gamma_\kappa(x)=\lim_{|x|\to 0}u(x)/\Gamma(x)=1$ for all $\kappa\in\R$, we have $v_\delta^-\leq u\leq v_\delta^+$ in a neighborhood of the origin. Finally, by the computations above we have $\De_p v_\delta^-\geq 0$ and $\De_p v_\delta^+\leq 0$.
Thus, by the
Comparison Principle for subsolutions and supersolutions of the p-Laplace equation (see, e.g., \cite[Theorem 3.5.1]{Puc_Ser_book} or \cite[Theorem 2.15]{LIN2017}),
we conclude 
\[
v_\delta^-\leq u\leq v_\delta^+
\]
on the whole punctured ball $B_{r_\ep}(o)\setminus \{o\}$.
Taking the limit as $\delta\to 0$, we conclude
\[
\Gamma_{a-\ep}(x)-\Gamma_{a-\ep}(r_\ep)+\inf_{\{|x|=r_\ep\}}u\leq u(x)\leq
\Gamma_{b+\ep}(x)-\Gamma_{b+\ep}(r_\ep)+\sup_{\{|x|=r_\ep\}}u\,.
\] 
Recalling the definition of $\Gamma_\kappa$, this inequality can be made more explicit, distinguishing three different cases. If $p\in(1,n)\setminus\{(n+2)/3\}$, we have obtained
\begin{equation}
\label{eq:error_bound_1}
\frac{a-\ep}{6(n+2-3p)}|x|^{-\frac{n+2-3p}{p-1}}+c_\ep\,\leq\, u(x)-\frac{p-1}{n-p}|x|^{-\frac{n-p}{p-1}}\,\leq\,\frac{b+\ep}{6(n+2-3p)}|x|^{-\frac{n+2-3p}{p-1}}+C_\ep\,,
\end{equation}
for all $x\in B_{r_\ep}(o)\setminus\{o\}$, where $c_\ep,C_\ep\in\R$.
If $p=(n+2)/3$, we have
\begin{equation}
\label{eq:error_bound_2}
\frac{a-\ep}{2(n-1)}\log\frac{1}{|x|}+c_\ep\,\leq\, u(x)-\frac{1}{2}|x|^{-2}\,\leq\,\frac{b+\ep}{2(n-1)}\log\frac{1}{|x|}+C_\ep\,.
\end{equation}
Finally, if $p=n$, it holds
\begin{equation}
\label{eq:error_bound_3}
-\frac{a-\ep}{12(n-1)}|x|^2+c_\ep\,\leq\, u(x)-\log\frac{1}{|x|}\,\leq\,-\frac{b+\ep}{12(n-1)}|x|^2+C_\ep\,.
\end{equation}
Estimates~\eqref{eq:error_bound_1},~\eqref{eq:error_bound_2} and~\eqref{eq:error_bound_3} in particular imply~\eqref{eq:assumption_error}.

\begin{remark}
\label{rmk:Einstein}
If $\Ric(o)$ is a multiple of the metric at the pole $o$, then choosing $\kappa=a=b$ yields, by the preceding computation,
$\De_p^{(g)}\Gamma_\kappa=\mathcal{O}(|x|^{3-n})$.
It is plausible that a sharper estimate for $u-\Gamma_\kappa$ could be obtained by adding a further correction term to $\Gamma_\kappa$, guided by the asymptotic expansion of the model solutions, estimating the corresponding $p$-Laplacian to the next order, and then proceeding as above.
\end{remark}

\section{Estimate on the H\"older norm of the error term}
\label{sec:estimates}

In this section we work in the Riemannian setting. The Euclidean setting will just be a special case.

Since $\{x^i\}$ are normal coordinates, we have $g_{ij}(x)=\delta_{ij}+\mathcal{O}_2(|x|^2)$. Let $s\in(0,1]$. We set
\[
g_s(x)\,=\,g_{ij}(sx)dx^i\otimes dx^j\,,\quad u_s(x)\,=\,s^{\frac{n-p}{p-1}}u(sx)\,,\quad e_s(x)\,=\,s^{\frac{n-p}{p-1}}e(sx)\,.
\]
Notice that this scaling does not affect the function $\Gamma$, since
\[
\Gamma(x)\,=\,s^{\frac{n-p}{p-1}}\Gamma(sx)\,.
\]

Let $R>0$. Notice that, if $\De_p^{(g)}u=0$ in $B_R\setminus\{o\}$, then
\[
\De^{(g_s)}_p u_s\,=\,0\ \hbox{ in }B_R\setminus\{o\}
\]
for all $s\in(0,1]$.

\subsection{\texorpdfstring{$C^1$}{C1}-estimates}
\label{sub:C1}

We recall from Section~\ref{sec:order_zero} that $\lim_{|x|\to 0}u(x)/\Gamma(x)=1$. 
As a consequence, given $\ep>0$, there exists $0<r_\ep<R$ such that
\begin{equation}
\label{eq:serrin-estimate}
(1-\ep)\Gamma(x)\,\leq\,u(x)\,\leq\,(1+\ep)\Gamma(x)\,,\quad\forall\ 0<|x|<r_\ep\,.
\end{equation}
Multiplying by $s^{\frac{n-p}{p-1}}$ and setting $x=sy$, we thus obtain
\[
(1-\ep)\Gamma(y)\,\leq\,u_s(y)\,\leq\,(1+\ep)\Gamma(y)\,,\quad\forall\ 0<|y|<r_\ep/s\,.
\]
As a consequence, there exists a constant $C=C(p,\ep)>0$ such that in the annulus $A=\{1/5\leq |x|\leq 4/5\}$ it holds
\[
\|u_s\|_{L^{\infty}(A)}\leq C
\]
for all $s\in(0,r_\ep]$. From~\cite{Tolksdorf} (see also~\cite[Theorem~B.3]{AgBoMa}) we conclude that there exist constants $0<\beta'<1,C'>0$ such that in $A'=\{1/4\leq |x|\leq 3/4\}$ it holds
\[
\|u_s\|_{C^{1,\beta'}(A')}\leq C'
\]
for all $s\in(0,r_\ep]$.
Given $0<\beta<\beta'$, it follows that, for any sequence $u_{s_k}$ with $s_k\to 0$ there exists a subsequence that converges to a function $\tilde u$ in $C^{1,\beta}(A')$.

Furthermore, for $s$ smaller than $r_\ep$, exploiting~\eqref{eq:serrin-estimate}, the error $e_s$ satisfies
\begin{equation}
\label{eq:es_Linfinity}
\|e_s\|_{L^{\infty}(A')}=s^{\frac{n-p}{p-1}}\|u(sx)-\Gamma(sx)\|_{L^\infty(A')}\,\leq\, s^{\frac{n-p}{p-1}}\sup_{y\in sA'}|\ep\Gamma(y)|\,=\,s^{\frac{n-p}{p-1}}\ep\Gamma({s}/{4})=C'' \ep\,,
\end{equation}
where $C''$ depends only on $p$. Since this holds for all $\ep>0$ and all $s<r_\ep$, it follows that $e_s(x)=u_s(x)-\Gamma(x)$ converges to zero in $L^\infty(A')$ as $s\to 0$, hence $\tilde u=\Gamma$.
In conclusion, we have proved that any sequence $u_{s_k}$ with $s_k\to 0$ admits a subsequence converging to $\Gamma$ in $C^{1,\beta}(A')$. As a consequence, it follows that $u_s$ converges to $\Gamma$ in  $C^{1,\beta}(A')$ as $s\to 0$ (in other words, $\|e_{s}\|_{C^{1,\beta}(A')}\to 0$ as $s\to 0$).

\subsection{\texorpdfstring{$C^2$}{C2}-estimates}

The aim of this subsection is to show that the error $e_s$ satisfies a uniformly elliptic equation and thus it is possible to apply Schauder estimates to control its $C^2$-norm. 

Given a point $x$ and a $1$-form $\xi$, let
\[
\mathcal{A}_s^i(x,\xi)\,=\,\left[\sqrt{g_s}\, g_s^{ij}\,|\xi|_{g_s}^{p-2}\xi_j\right](x)\,.
\]
We can write
\begin{align*}
0\,&=\,\De_p^{(g_s)}u_s\,=\,\De_p^{(g_s)}(\Gamma+e_s)-\De_p^{(g_s)}\Gamma+\De_p^{(g_s)}\Gamma
\\
&=\,\frac{1}{\sqrt{g_s}}\,\frac{\pa}{\pa x^i}\mathcal{A}_s^i(x,d\Gamma+de_s)-\frac{1}{\sqrt{g_s}}\,\frac{\pa}{\pa x^i}\mathcal{A}_s^i(x,d\Gamma)+\De_p^{(g_s)}\Gamma\,,
\end{align*}
that is
\begin{equation}
\label{eq:check_elliptic}
\frac{\pa}{\pa x^i}\left[\mathcal{A}_s^i(x,d\Gamma+de_s)-\mathcal{A}_s^i(x,d\Gamma)\right]\,=\,-\left[\sqrt{g_s}\,\De_p^{(g_s)}\Gamma\right](x)\,. 
\end{equation}
We now proceed to study the left-hand side of~\eqref{eq:check_elliptic}, showing that it is a uniformly elliptic operator applied to $e_s$. To this end, observe that
\[
\mathcal{A}_s^i(x,d\Gamma+de_s)-\mathcal{A}_s^i(x,d\Gamma)\,=\,\int_0^1\frac{d}{d t}\mathcal{A}_s^i(x,d\Gamma+t de_s)dt\,=\,\int_0^1\frac{\pa\mathcal{A}_s^i}{\pa \xi^j}(x,d\Gamma+tde_s)\frac{\pa e_s}{\pa x^j}(x)dt\,.
\]
Thus, we can rewrite~\eqref{eq:check_elliptic} as
\begin{equation}
\label{eq:elliptic}
\frac{\pa}{\pa x^i}\left(a_s^{ij}(x)\frac{\pa e_s}{\pa x^j}(x)\right)\,=\,-[\sqrt{g_s}\,\De_p^{(g_s)}\Gamma](x)\,,
\end{equation}
with
\[
a_s^{ij}(x)\,=\,\int_0^1\frac{\pa \mathcal{A}_s^i}{\pa \xi^j}(x,d\Gamma+tde_s)dt\,.
\]
A simple computation shows that
\[
\frac{\pa\mathcal{A}_s^i}{\pa \xi^j}(x,\xi)\,=\,
\left\{\sqrt{g_s}|\xi|_{g_s}^{-(4-p)}\,\left[|\xi|_{g_s}^2 g^{ij}_s-(2-p)g_s^{ik}g_s^{jm}\xi_k\xi_m\right]\right\}(x)
\]
hence in particular, for every $1$-form $\eta$ we compute
\begin{align*}
a_s^{ij}(x)\eta_i\eta_j\,&=\,\int_0^1\sqrt{g_s}|d\Gamma+tde_s|_{g_s}^{-(4-p)}\,\left[|d\Gamma+tde_s|_{g_s}^2 |\eta|^2_{g_s}-(2-p)\langle d\Gamma+tde_s,\eta\rangle_{g_s}^2\right]dt
\\
&\geq\,\min\{p-1,1\}\sqrt{g_s} |\eta|^2_{g_s}\int_0^1|d\Gamma+tde_s|_{g_s}^{-(2-p)}\,dt\,,
\end{align*}
where in the latter inequality we have used Cauchy--Schwartz if $p<2$, whereas we have just neglected the last term if $p\geq 2$. 
We now recall from Subsection~\ref{sub:C1} that $\|e_{s}\|_{C^1(A')}\to 0$ as $s\to 0$. Furthermore, for small values of $s$, the metric $g_s$ is close to the Euclidean metric and the derivatives of $\Gamma$ are bounded in $A'$, hence there exist positive constants $c_1,c_2$ and $0<\bar s<1$ such that, for all $s\in(0,\bar s]$, it holds
\[
c_1\leq\frac{|d\Gamma|_{g_{\R^n}}}{2}\leq |d\Gamma+tde_s|_{g_s}\leq 2|d\Gamma|_{g_{\R^n}}\leq c_2
\]
at all points of $A'$.
Thus, we conclude that there exists a constant $\lambda$, depending on $p$ but independent of $s$, such that $a_s^{ij}\eta_i\eta_j\geq\lambda|\eta|_{g_s}^2$ for all $s\in(0,\bar s]$. With completely analogous estimates we find also the opposite inequality $a_s^{ij}\eta_i\eta_j\leq\Lambda|\eta|_{g_s}^2$.
In conclusion, the left-hand side of~\eqref{eq:elliptic} is a uniformly elliptic operator in $A'$ applied to the error term $e_s$, for all $s\in(0,\bar s]$.

Finally, observe that the $C^{1,\beta}$-norm of $e_s$ is bounded by a constant independent of $s$ and that $\Gamma$ and its derivatives are also bounded in $A'$. As a consequence of this and the explicit expression of $a_s^{ij}$, it is easily seen that the $C^{0,\beta}$-norm of the coefficients $a_s^{ij}$ is bounded by a constant independent of $s$.

Therefore, by the Schauder interior estimates~\cite[Corollary 6.3]{GiTr}, it holds
\begin{equation}
\label{eq:schauder}
\left\|e_s\right\|_{C^{2,\beta}(A'')}\,\leq\,C'''\left(\left\|e_s\right\|_{L^{\infty}(A')}+\|\sqrt{g_s}\,\De_p^{(g_s)}\Gamma\|_{C^{0,\beta}(A')}\right)\,,
\end{equation}
where $A''=\{1/3\leq|x|\leq 2/3\}$, $s\in(0,\bar s]$ and $C'''$ is a constant that depends on $p$ but not on $s$.

\section{Proof of the main theorems}

\subsection{Proof of Theorem~\ref{thm:Riemannian_1}}

To exploit~\eqref{eq:schauder}, we first need to estimate $\|\sqrt{g_s}\,\De_p^{(g_s)}\Gamma\|_{C^{0,\beta}(A')}$. To this end, we compute
\begin{multline*}
\sqrt{g_s}\,\De_p^{(g_s)}\Gamma\,=\,\frac{\pa}{\pa x^i}\mathcal{A}_s^i(x,d\Gamma)
\\
=\,\frac{\pa\sqrt{g_s}}{\pa x^i} g_s^{ij}|d\Gamma|_{g_s}^{p-2}\frac{\pa\Gamma}{\pa x^j}+
\sqrt{g_s} \frac{\pa g_s^{ij}}{\pa x^i}|d\Gamma|_{g_s}^{p-2}\frac{\pa\Gamma}{\pa x^j}+
\sqrt{g_s} g_s^{ij}\frac{\pa |d\Gamma|_{g_s}^{p-2}}{\pa x^i}\frac{\pa\Gamma}{\pa x^j}+
\sqrt{g_s} g_s^{ij}|d\Gamma|_{g_s}^{p-2}\frac{\pa^2\Gamma}{\pa x^i\pa x^j}\,.
\end{multline*}
Notice that $(g_s)_{ij}=\delta_{ij}+\mathcal{O}(s^2)$ and $\pa (g_s)_{ij}/\pa x^k=\mathcal{O}(s^2)$. Furthermore, $\Gamma$ and its derivatives are uniformly bounded in the annulus $A'$. As a consequence:
\[
\sqrt{g_s}\,\De_p^{(g_s)}\Gamma\,=\, \delta^{ij}\frac{\pa |d\Gamma|_{g_s}^{p-2}}{\pa x^i}\frac{\pa\Gamma}{\pa x^j}+
\delta^{ij}|d\Gamma|_{g_s}^{p-2}\frac{\pa^2\Gamma}{\pa x^i\pa x^j}+\mathcal{O}(s^2)=\De_p^{(g_{\R^n})}\Gamma+\mathcal{O}(s^2)\,=\,\mathcal{O}(s^2)\,.
\]
Proceeding in the same way, we can also estimate the derivative:
\[
\frac{\pa}{\pa x^i}\left[\sqrt{g_s}\,\De_p^{(g_s)}\Gamma\right]\,=\,\frac{\pa}{\pa x^i}\De_p^{(g_{\R^n})}\Gamma+\mathcal{O}(s^2)\,=\,\mathcal{O}(s^2)\,.
\]
In particular, we have shown that there exists a constant $c$ such that
\begin{equation}
\label{eq:estimate_Deltap-Gamma}
\|\sqrt{g_s}\De_p^{(g_s)}\Gamma\|_{C^{0,\beta}(A')}\leq\|\sqrt{g_s}\De_p^{(g_s)}\Gamma\|_{C^1(A')}\leq cs^2\,.
\end{equation}

Recalling~\eqref{eq:assumption_error}, we can also refine~\eqref{eq:es_Linfinity}:
\begin{equation*}
\|e_s\|_{L^{\infty}(A')}=s^{\frac{n-p}{p-1}}\|e(sx)\|_{L^\infty(A')}\leq s^{\frac{n-p}{p-1}}\sup_{y\in sA'}|e(y)|\leq 
\begin{dcases}
C'' s^{\frac{n-p}{p-1}}& \hbox{if }\frac{n+2}{3}<p<n
\\
C'' s^2 \log\frac{1}{s}  & \hbox{if } p=\frac{n+2}{3}
\\
C'' s^2  & \hbox{if }1<p<\frac{n+2}{3}
\end{dcases}
\end{equation*}

We now study the three cases separately.
If $p<\frac{n+2}{3}$, the above estimate, together with~\eqref{eq:schauder} and~\eqref{eq:estimate_Deltap-Gamma}, give
\[
\left\|e_s\right\|_{C^2(A'')}\,\leq\,K s^2\,,
\]
or equivalently:
\[
s^{\frac{n-p}{p-1}}|e(sx)|
+s^{\frac{n-p}{p-1}+1}\sum_{i=1}^n\left|\frac{\pa e}{\pa x^i}(sx)\right|+s^{\frac{n-p}{p-1}+2}\sum_{i,j=1}^n\left|\frac{\pa^2 e}{\pa x^i\pa x^j}(sx)\right|\,\leq\,Ks^2
\]
for all $x\in A''$. Notice that, for all $x\in A''$ it holds $1/3\leq|x|\leq 2/3$. In particular $|x|\leq 1$, and so $s\geq |sx|$ and $s^{2-\frac{n-p}{p-1}}\leq |sx|^{2-\frac{n-p}{p-1}}$. Setting $y=sx$, the previous inequality gives
\begin{equation}
\label{eq:thesis1}
|e(y)|
+|y|\sum_{i=1}^n\left|\frac{\pa e}{\pa x^i}(y)\right|+|y|^2\sum_{i,j=1}^n\left|\frac{\pa^2 e}{\pa x^i\pa x^j}(y)\right|\,\leq\,K|y|^{2-\frac{n-p}{p-1}}\,.
\end{equation}
for all $y$ with $s/3\leq |y|\leq 2s/3$. Since this is true for all $s\in(0,\bar s]$, then it is true for all $y\in B(o,2\bar s/3)\setminus\{o\}$.

\smallskip

If instead $p> \frac{n+2}{3}$, then $\frac{n-p}{p-1}< 2$ and so, employing again~\eqref{eq:schauder} and~\eqref{eq:estimate_Deltap-Gamma}, we have
\[
\left\|e_s\right\|_{C^2(A'')}\,\leq\,K s^{\frac{n-p}{p-1}}\,,
\]
or equivalently:
\[
s^{\frac{n-p}{p-1}}|e(sx)|
+s^{\frac{n-p}{p-1}+1}\sum_{i=1}^n\left|\frac{\pa e}{\pa x^i}(sx)\right|+s^{\frac{n-p}{p-1}+2}\sum_{i,j=1}^n\left|\frac{\pa^2 e}{\pa x^i\pa x^j}(sx)\right|\,\leq\,Ks^{\frac{n-p}{p-1}}
\]
for all $x\in A''$. Notice that, for all $x\in A''$ it holds $1/3\leq|x|\leq 2/3$. In particular $s\geq |sx|$. Setting $y=sx$, the previous inequality gives
\begin{equation}
\label{eq:thesis2}
|e(y)|
+|y|\sum_{i=1}^n\left|\frac{\pa e}{\pa x^i}(y)\right|+|y|^2\sum_{i,j=1}^n\left|\frac{\pa^2 e}{\pa x^i\pa x^j}(y)\right|\,\leq\,K
\end{equation}
for all $y$ with $s/3\leq |y|\leq 2s/3$. Since this is true for all $s\in(0,\bar s]$, then it is true for all $y\in B(o,2\bar s/3)\setminus\{o\}$.

\smallskip

Finally, if $p=\frac{n+2}{3}$, then from~\eqref{eq:schauder} and~\eqref{eq:estimate_Deltap-Gamma} we pbtain
\[
\left\|e_s\right\|_{C^2(A'')}\,\leq\,K s^2\log(1/s)\,,
\]
or equivalently:
\[
s^2|e(sx)|
+s^3\sum_{i=1}^n\left|\frac{\pa e}{\pa x^i}(sx)\right|+s^4\sum_{i,j=1}^n\left|\frac{\pa^2 e}{\pa x^i\pa x^j}(sx)\right|\,\leq\,Ks^2\log(1/s)
\]
for all $x\in A''$. Notice that, for all $x\in A''$ it holds $1/3\leq|x|\leq 2/3$. In particular $s\geq |sx|$. Setting $y=sx$, the previous inequality gives
\begin{equation}
\label{eq:thesis3}
|e(y)|
+|y|\sum_{i=1}^n\left|\frac{\pa e}{\pa x^i}(y)\right|+|y|^2\sum_{i,j=1}^n\left|\frac{\pa^2 e}{\pa x^i\pa x^j}(y)\right|\,\leq\,K\log(1/|y|)
\end{equation}
for all $y$ with $s/3\leq |y|\leq 2s/3$. Since this is true for all $s\in(0,\bar s]$, then it is true for all $y\in B(o,2\bar s/3)\setminus\{o\}$. 
Estimates~\eqref{eq:thesis1},~\eqref{eq:thesis2} and~\eqref{eq:thesis3} are precisely the thesis of Theorem~\ref{thm:Riemannian_1}.

\subsection{Proof of Theorem~\ref{thm:Euclidean_1}}

In the setting of Theorem~\ref{thm:Euclidean_1}, the metric is Euclidean, that is, $g_{ij}=\delta_{ij}$, without any error term (in particular $g_s=g$). 
Furthermore, from~\cite{AnEs} we know that in the Euclidean setting the error term $e$ is $C^{0,\alpha}$. In particular, we have
\begin{equation}
\label{eq:e_Linfty}
\left|e_s(x)-s^{\frac{n-p}{p-1}}e(o)\right|=s^{\frac{n-p}{p-1}}|e(sx)-e(o)|\leq
\widetilde{C} s^{\frac{n-p}{p-1}}|sx|^\alpha
\leq
\widetilde{C} s^{\frac{n-p}{p-1}+\alpha}
\end{equation}
in $A'$.

The estimates in Section~\ref{sec:estimates} still apply, but it is convenient to rewrite~\eqref{eq:elliptic} in the following clearly equivalent form:
\[
\frac{\pa}{\pa x^i}\left[a_s^{ij}(x)\frac{\pa}{\pa x^j}\left(e_s(x)-s^{\frac{3-p}{p-1}}e(o)\right)\right]\,=\,-\De_p^{(g_{\R^n})}\Gamma\,.
\]
Furthermore, notice that, in the Euclidean setting, $\Gamma$ is $p$-harmonic, hence $\De_p^{(g_{\R^n})}\Gamma=0$.
This way, applying Schauder estimates as in~\eqref{eq:schauder}, we get
\[
\left\|e_s-s^{\frac{n-p}{p-1}}e(o)\right\|_{C^2(A'')}\,\leq\,C'''\left\|e_s-s^{\frac{n-p}{p-1}}e(o)\right\|_{L^{\infty}(A')}\leq\,Ks^{\frac{n-p}{p-1}+\alpha}\,,
\]
where in the last inequality we have used~\eqref{eq:e_Linfty}.
Equivalently, it holds
\[
s^{\frac{n-p}{p-1}}|e(sx)-e(o)|
+s^{\frac{n-p}{p-1}+1}\sum_{i=1}^n\left|\frac{\pa e}{\pa x^i}(sx)\right|+s^{\frac{n-p}{p-1}+2}\sum_{i,j=1}^n\left|\frac{\pa^2 e}{\pa x^i\pa x^j}(sx)\right|\,\leq\,Ks^{\frac{n-p}{p-1}+\alpha}
\]
for all $x\in A''$. Notice that, for all $x\in A''$ it holds $1/3\leq|x|\leq 2/3$. In particular $s\geq |sx|$ and $s^\alpha\leq 3|sx|^\alpha$. Setting $y=sx$, the previous inequality gives
\[
|e(y)-e(o)|
+|y|\sum_{i=1}^n\left|\frac{\pa e}{\pa x^i}(y)\right|+|y|^2\sum_{i,j=1}^n\left|\frac{\pa^2 e}{\pa x^i\pa x^j}(y)\right|\,\leq\,3K|y|^\alpha\,.
\]
for all $y$ with $s/3\leq |y|\leq 2s/3$. Since this is true for all $s\in(0,\bar s]$, then it is true for all $y\in B(o,2\bar s/3)\setminus\{o\}$.
This is precisely our thesis.

\appendix

\section{The model solutions and their expansion near the pole}
\label{appendix}

Let $(M_\kappa,g_\kappa)$ be a space form with sectional curvatures equal to $\frac{\kappa}{n-1}$ (so that in particular $\Ric_{g_\kappa}=\kappa g_\kappa$) and let $x_1,\dots,x_n$ be normal coordinates centered at $o\in M_\kappa$. The metric $g_\kappa$ can be written as
\[
g_\kappa\,=\,d|x|\otimes d|x|+\operatorname{sn}^2_\kappa |x|\, g_{\mathbb{S}^{n-1}}\,,
\]
where 
\[
\operatorname{sn}_\kappa\rho\,=\,
\begin{cases}
\sqrt{\frac{n-1}{\kappa}} \sin\left(\sqrt{\frac{\kappa}{n-1}}\,\rho\right) & \hbox{if }\kappa>0\,,
\\
\rho & \hbox{if }\kappa=0\,,
\\
\sqrt{\frac{n-1}{|\kappa|}} \sinh\left(\sqrt{\frac{|\kappa|}{n-1}}\,\rho\right) & \hbox{if }\kappa<0\,.
\end{cases}
\]
Given $p\in(1,n]$, it is then easy to compute that the rotationally symmetric $p$-harmonic functions are given by
\[
u_\kappa(x)\,=\,\int_{|x|}^{c}\left(\operatorname{sn}_\kappa\rho\right)^{-\frac{n-1}{p-1}} d\rho\,,
\]
with $c>0$. Furthermore, we have the following Taylor expansion for $\rho\to 0$: 
\begin{equation}
\label{eq:expansion_sn}
\left(\operatorname{sn}_\kappa\rho\right)^{-\frac{n-1}{p-1}}\,=\,\rho^{-\frac{n-1}{p-1}}+\frac{\kappa}{6(p-1)}\rho^{2-\frac{n-1}{p-1}}+\mathcal{O}(\rho^{4-\frac{n-1}{p-1}})\,.
\end{equation}
Integrating between $|x|$ and $c$, if $p\neq(n+2)/3$ and $p\neq (n+4)/5$, it is then easy to conclude
\[
u_\kappa(x)\,=\,
\frac{p-1}{n-p}|x|^{-\frac{n-p}{p-1}}+\frac{\kappa}{6(n+2-3p)}|x|^{-\frac{n+2-3p}{p-1}}+C+\mathcal{O}(|x|^{-\frac{n+4-5p}{p-1}})\,,
\]
where the contant $C$ depends on $c$. Up to adding a constant to $u_\kappa$, we can always set $C$ to be equal to zero.

In the case $p=(n+4)/5$, the error term in~\eqref{eq:expansion_sn} is $\mathcal{O}(1/\rho)$, hence when we integrate it we obtain a constant plus an error term that grows as $\log|x|$. Thus, we obtain the estimate
\[
u_\kappa(x)\,=\,
\frac{1}{4}|x|^{-4}+\frac{5\kappa}{12(n-1)}|x|^{-2}+C+\mathcal{O}(\log(1/|x|))\,.
\]
In the case $p=(n+2)/3$, the second term in~\eqref{eq:expansion_sn} grows as $1/\rho$, hence we compute
\[
u_\kappa(x)\,=\,
\frac{1}{2}|x|^{-2}+\frac{\kappa}{2(n-1)}\log\frac{1}{|x|}+C+\mathcal{O}(|x|^2)\,.
\]
Finally, in the case $p=n$, the first term in~\eqref{eq:expansion_sn} is equal to $1/\rho$, therefore
\[
u_\kappa(x)\,=\,
\log\frac{1}{|x|}-\frac{\kappa}{12(n-1)}|x|^2+C+\mathcal{O}(|x|^4)\,.
\]

\begin{ackn}
The authors would like to thank L.~Mazzieri for useful discussions during the preparation of this work. The authors are members of the Gruppo Nazionale per l’Analisi Matematica, la Probabilità e le loro Applicazioni (GNAMPA), which is part of the Istituto Nazionale di Alta Matematica (INdAM). 
S.~B. is partially funded by the GNAMPA project ``Analysis of Non-smooth Geometric Evolutions''.
\end{ackn}

\bibliographystyle{plain}
\bibliography{biblio}

\end{document}